\documentclass[10pt]{amsart}
\usepackage{amsmath,amssymb,amsfonts, amscd,enumerate}
\usepackage{amsbsy}
\usepackage{amsthm}
\usepackage{amstext}
\usepackage{amsopn}
\usepackage{hyperref}
\usepackage{mathrsfs}
\usepackage{graphicx}
\usepackage{latexsym}
\usepackage[T1]{fontenc}
\usepackage{color}
\usepackage{comment}
\usepackage{pgflibraryarrows}
\usepackage{pgflibrarysnakes}
\usepackage{tikz, ifthen, multirow}

\newcommand{\Hmm}[1]{\leavevmode{\marginpar{\tiny%
$\hbox to 0mm{\hspace*{-0.5mm}$\leftarrow$\hss}%
\vcenter{\vrule depth 0.1mm height 0.1mm width \the\marginparwidth}%
\hbox to 0mm{\hss$\rightarrow$\hspace*{-0.5mm}}$\\\relax\raggedright #1}}}

\newcommand{\Dr}{\mathscr{D}}

\newcommand{\Tr}{\mathscr{T}}

\newtheorem{theorem}{Theorem}

\newtheorem{thm}{Theorem}[section]

\newtheorem{lemma}[thm]{Lemma}
\newtheorem{pro}[thm]{Proposition}

\theoremstyle{definition}
\newtheorem*{defi}{Definition}

\newtheorem{rem}[thm]{Remark}
\numberwithin{equation}{section}

\newcommand{\Z}{{\mathbb Z}}

\begin{document}
\title[]{On the uniqueness of solutions for the Fubuki game}

\author[S. Gol\'enia]{Sylvain Gol\'enia} \address{Institut de
  Math\'ematiques de Bordeaux Universit\'e Bordeaux 1
351, cours de la Lib\'eration
F-33405 Talence cedex, France} \email{sylvain.golenia@math.u-bordeaux1.fr}


\date{\today}

\maketitle
\begin{abstract} \noindent
We discuss the problem of uniqueness of solutions for the
Fubuki game when the diagonal is prescribed.
\end{abstract}

\section{Introduction}
The Fubuki game is a generalization of the magical square. The rules
are simple: Fill a $3$ by $3$ grid with the numbers $1$ to $9$ so that
each line and column adds up to a given sum. For instance,
\begin{align*}
\begin{tabular}{|c|c|c|r|}
\hline
{\bf 1} & & & =10
\\
\hline
& {\bf 2 } & &=15
\\
\hline
& &  {\bf 3} & =20
\\
\hline
= 16 & =15 & =14
\\
\cline{1-3}
\end{tabular} 
\end{align*}
leads to the two solutions:
\begin{align*}
\begin{tabular}{|c|c|c|r|}
\hline
{\bf 1} & 4&5 & =10
\\
\hline
7 & {\bf 2 } & 6&=15
\\
\hline
8 & 9 &  {\bf 3} & =20
\\
\hline
= 16 & =15 & =14
\\
\cline{1-3}
\end{tabular}
\quad\quad
 \begin{tabular}{|c|c|c|r|}
\hline
{\bf 1} & 5 &4  & =10
\\
\hline
6 & {\bf 2 } & 7 &=15
\\
\hline
9 &8 &  {\bf 3} & =20
\\
\hline
= 16 & =15 & =14
\\
\cline{1-3}
\end{tabular}
\end{align*}
whereas 
\begin{align*}
\begin{tabular}{|c|c|c|r|}
\hline
{\bf 1} & & & =11
\\
\hline
& {\bf 2 } & &=14
\\
\hline
& &  {\bf 3} & =20
\\
\hline
= 14 & =15 & =16
\\
\cline{1-3}
\end{tabular} 
\end{align*}
leads to the unique solution:
\begin{align*}
\begin{tabular}{|c|c|c|r|}
\hline
{\bf 1} & 4&6 & =11
\\
\hline
5 & {\bf 2 } & 7&=14
\\
\hline
8 & 9 &  {\bf 3} & =20
\\
\hline
= 14 & =15 & =16
\\
\cline{1-3}
\end{tabular}
\end{align*}
 
In these two examples we have prescribed the diagonal. It is
not necessary in the general case but it will be the setting of our
study. It is the first interesting case since, as a linear system,
i.e, solving it with real numbers, the affine subspace of solutions is
of dimension $1$, when the solution exists.

The first problem is to ensure the uniqueness of the
solution. This can be done solely by choosing the diagonal.

\begin{theorem}\label{t:unicity}
a) Assume that the diagonal belongs to the following set:
\begin{align*}
\Dr:=\big\{&
\{1,3,4\}, \{1,3,6\}, \{1,4,6\}, \{1,5,7\}, \{1,5,8\}, \{1,7,8\}, \{2,3,4\},
\{2,3,5\}, 
\\
&\{2,3,6\}, \{2,3,8\}, \{2,3,9\}, \{2,4,6\}, \{2,4,7\}, \{2,4,8\},
\{2,5,6\}, \{2,5,9\}, 
\\
&\{2,6,7\}, \{2,6,8\}, \{2,7,8\}, \{3,4,6\}, \{3,4,8\},
\{3,5,6\}, \{3,5,7\}, \{3,5,9\}, 
\\
&\{3,6,8\}, \{4,5,7\}, \{4,5,8\},
\{4,6,7\}, \{4,6,8\}, \{4,6,9\}, \{4,7,8\}, \{4,7,9\},
\\
& \{5,7,8\}, \{6,7,8\},
\{6,7,9\}\big\}.
\end{align*}
Then if the solution of the Fubuki game with prescribed diagonal exists,
it is unique. 
\\
b) If the diagonal does not belong to $\Dr$, there are at most $2$
solutions. 
\end{theorem} 
We stress that the result is independent of the choice of the
sum of the lines and that of the column.

Next we count the number of grids that lead to a unique solution. 

\begin{theorem}\label{t:main}
Among the Fubuki grids with prescribed diagonal which have a solution,
there are $351432$ of them which have a unique solution. 
\end{theorem} 
Recalling that there are $362880$ possible grids we have therefore
around $96\%$ chance to have a unique solution. The novelty in this
note is not really about the statement of the two theorems, as it
could have been guessed and done pretty easily with Matlab for instance, but
relies in the fact that it can be humanly proven. 

To conclude we mention that brute force computations give  that there are:
$281304(\simeq 77 \%)$ unique solutions if we prescribe the two first
squares on the diagonal, $163387(\simeq 45 \%)$ ones if we prescribe
the upper left square, and $46147 (\simeq 12 \%)$ ones if we prescribe
no square at all.

\section{Proofs}
Suppose that we have a solution of a Fubuki
problem with prescribed diagonal given by
\begin{align}\label{e:sol}
S_0:=\begin{tabular}{|c|c|c|}
\hline
${\bf s_{1,1}}$& $s_{1,2}$ & $s_{1,3}$
\\
\hline
$s_{2,1}$& ${\bf s_{2,2}}$ & $s_{2,3}$
\\
\hline
$s_{3,1}$& $s_{3,2}$ & ${\bf s_{3,3}}$
\\
\hline
\end{tabular}
\end{align}
We stress that we have prescribed the diagonal by using a bold font.

If we have a second solution then
there is $a\in \Z\setminus\{0\}$ such that $s_{1,2}+a$ is at the place of
the  $s_{1,2}$. Since the sum of the elements of the first line is the
same for the two grids, we see that $s_{1,3}$ is hence replaced by
$s_{1,3}-a$. We fill the rest of the grid in the same manner, we
obtain:

\begin{align}\label{e:sol2}
S_a:=\begin{tabular}{|c|c|c|}
\hline
${\bf s_{1,1}}$& $s_{1,2}+a$ & $s_{1,3}-a$
\\
\hline
$s_{2,1}-a$& ${\bf s_{2,2}}$ & $s_{2,3}+a$
\\
\hline
$s_{3,1}+a$& $s_{3,2}-a$ & ${\bf s_{3,3}}$
\\
\hline
\end{tabular}
\end{align}
Note that $S_a$ is an other solution if and only if
\begin{align}
\nonumber
\{s_{1,2}+a, s_{1,3}-a,s_{2,1}-a, s_{2,3}+a, s_{3,1}+a, s_{3,2}-a
\}&
\\
\label{e:firsta}
&\hspace*{-3cm}=\{s_{1,2}, s_{1,3},s_{2,1}, s_{2,3}, s_{3,1}, s_{3,2}\}:=X.
\end{align}

We give a first remark on the structure of the solutions.

\begin{rem}\label{r:sign}
By considering the minimal element of the set $X$ we see that: if
there are $a,b\in \Z\setminus\{0\}$ such that $S_a$ and $S_b$ exist,
then $a$ and $b$  are of the same sign.
\end{rem}

\subsection{Triplet structure}

We turn to the structure of the set $X$ with respect to $a$. 

\begin{lemma}\label{l:y+a}
Given $a\in \Z\setminus \{0\}$ such that \eqref{e:firsta} holds true,
then there is a unique triplet 
$y_1, y_2, y_3 \in X$ such that $y_1<y_2<y_3$ and
\[\{y_1,  y_2,  y_3\}+ \{0, a\}:=\{y_1,
y_1+a, y_2, y_2+a, y_3, y_3+a\}= X.\]  
Here the plus between ensembles denotes the Minkowski sum.
\end{lemma}
\proof
To fix ideas say that $a>0$. Set $y_1:= \min (X)$.  Therefore by
\eqref{e:firsta} $y_1+a$ or $y_1-a$ belongs to $X$. By
minimality, we have that $y_1+a\in X$. To complete the proof, take
$y_2:= \min (X  \setminus \{y_1, y_1+a\})$ and $y_3:=
\min(X\setminus\{y_1, y_1+a, y_2, y_2+a\})$. The uniqueness follows
by construction. \qed

\begin{defi}
When it exists, we call $\Tr_a:=\{y_1, y_2, y_3\}$ the \emph{triplet}
associated to $a$ and to $S_0$.
\end{defi} 
 
\begin{rem}\label{r:reduc}
Note that $\Tr_a$ exists if and only if $\Tr_{-a}$ exists (set
$y_i(a):=y_i(-a)+a$). Therefore to prove the existence of $\Tr_a$ it
is enough to do it for $a>0$. However, recalling Remark \ref{r:sign},
$S_a$ cannot be a solution if $S_{-a}$ is one.
\end{rem} 

Next we notice the incompatibility of the two different triplets.

\begin{lemma}\label{l:diff}
Given $a$ and $b$ in $\Z\setminus \{0\}$ such that \eqref{e:sol2} and
\eqref{e:firsta} hold true and that $a\neq b$. Then $\Tr_a\neq \Tr_b$. 
\end{lemma} 
\proof First note that by Remark \ref{r:sign} that $a$ and $b$ are of
the same sign. To fix ideas say that they are both positive.
 Suppose
that $\Tr_a=\Tr_b$. By maximality $y_3(a)+a= y_3(b)+b$. We infer
that $a=b$. Contradiction. \qed

Finally we rephrase the existence of the solution \eqref{e:sol2} with the help
of the triplet $\Tr_a$ into a key-stone Proposition.

\begin{pro}\label{p:suite}
Let $S_0$ be the solution of \eqref{e:sol} and $a\in \Z$.
\\
a) If the triplet $\Tr_a$ associated to  $S_0$ and to $a\in
\Z\setminus\{0\}$ does not exists, then $S_a$ does not exist. 
\\
b) Assume that there exists the triplet $\Tr_a$ associated to  $S_0$
and to $a\in \Z\setminus\{0\}$. Then the solution $S_a$ exists if and
only if:  
\begin{align}\label{e:a>0}
\Tr_a= \{s_{1,2}, s_{2,3}, s_{3,1}\}.
\end{align}
\end{pro} 
\proof The point a) is clear from \eqref{e:firsta} and Lemma \ref{l:y+a}. 

We
turn to b). To fix idea say that $a>0$. We denote by $S^+:=\{s_{1,2},
s_{2,3}, s_{3,1}\}$ and by $S^-:=X\setminus S^+$. If $\Tr_a=S^+$ by
Lemma \ref{l:y+a} we see that $\Tr_a+\{a\}= S^-$. Consequently
\eqref{e:sol2} and \eqref{e:firsta} are satisfied: $S_a$ exists. 

Assume now that $S_a$ exists: \eqref{e:sol2} and \eqref{e:firsta} are
satisfied. If $y_3+a\in S^+$ then we obtain a contradiction by
maximality on $X$ and by \eqref{e:firsta}. Therefore $y_3+a\in
S^-$. Next since $y_2<y_3$, \eqref{e:firsta} ensures that $y_3\in
S^+$. We repeat the proof for $y_2$ and then for $y_1$ to conclude that
$\Tr_a= S^+$. 
\qed

\subsection{List of possibilities}\label{s:list}
As seen in Proposition \ref{p:suite}, the number of solutions of a
grid is linked with the existence of a Triplet. We classify
them. Given 
$x_1,x_2,x_3\in \{1,\ldots, 9\}$, in the first column, we give the
possible $c>0$, in the second column, such that there exist $y_1,
y_2$, and $y_3$ in $\{1,\ldots, 9\}\setminus\{x_1,x_2,x_3\}$ satisfying 
\begin{align}\label{e:a+}
\{y_1, y_1+c, y_2, y_2+c, y_3, y_3+c\}= \{1,\ldots, 9\}\setminus\{x_1,x_2,x_3\}.
\end{align}
This gives:
\begin{align*}
\begin{tabular}{|c|c|c|c|c|c|c|c|c|c|c|}
\cline{1-2} \cline{4-5} \cline{7-8} \cline{10-11}
                      \{1, 2, 3\} & \{1, 3\} & &
                         \{1, 2, 4\} & \{2\} & &
                         \{1, 2, 5\}& \{1\} & &
                         \{1, 2, 6\}& \{4\}
\\
\cline{1-2} \cline{4-5} \cline{7-8} \cline{10-11}
                                    \{1, 2, 7\}& \{1\}& &
                         \{1, 2, 8\}& \{2\}& &
                       \{1, 2, 9\}& \{1, 3\}& &
                         \{1, 3, 4\}& $\emptyset$
\\
\cline{1-2} \cline{4-5} \cline{7-8} \cline{10-11}
                         \{1, 3, 5\}& \{2\}&&
                         \{1, 3, 6\}& $\emptyset$&&
                         \{1, 3, 7\}& \{4\}&&
                         \{1, 3, 8\}& \{3\}
\\
\cline{1-2} \cline{4-5} \cline{7-8} \cline{10-11}
                         \{1, 3, 9\}& \{2\}&&
                         \{1, 4, 5\}& \{1\}&&
                         \{1, 4, 6\}& $\emptyset$&&
                         \{1, 4, 7\}& \{1\}
\\
\cline{1-2} \cline{4-5} \cline{7-8} \cline{10-11}
                         \{1, 4, 8\}& \{4\}&&
                         \{1, 4, 9\}& \{1\}&&
                         \{1, 5, 6\}& \{5\}&&
                         \{1, 5, 7\}& $\emptyset$
\\
\cline{1-2} \cline{4-5} \cline{7-8} \cline{10-11}
                         \{1, 5, 8\}& $\emptyset$&&
                         \{1, 5, 9\}& \{4\}&&
                         \{1, 6, 7\}& \{1\}&&
                         \{1, 6, 8\}& \{2\}
\\
\cline{1-2} \cline{4-5} \cline{7-8} \cline{10-11}
                         \{1, 6, 9\}& \{1\}&&
                         \{1, 7, 8\}& $\emptyset$&&
                         \{1, 7, 9\}& \{2\}&&
                       \{1, 8, 9\}& \{1, 3\}
\\
\cline{1-2} \cline{4-5} \cline{7-8} \cline{10-11}
                         \{2, 3, 4\}& $\emptyset$&&
                         \{2, 3, 5\}& $\emptyset$&&
                         \{2, 3, 6\}& $\emptyset$&&
                         \{2, 3, 7\}& \{3\}
\\
\cline{1-2} \cline{4-5} \cline{7-8} \cline{10-11}
                         \{2, 3, 8\}& $\emptyset$&&
                         \{2, 3, 9\}& $\emptyset$&&
                         \{2, 4, 5\}& \{2\}&&
                         \{2, 4, 6\}& $\emptyset$
\\
\cline{1-2} \cline{4-5} \cline{7-8} \cline{10-11}
                         \{2, 4, 7\}& $\emptyset$&&
                         \{2, 4, 8\}& $\emptyset$&&
                         \{2, 4, 9\}& \{2\}&&
                         \{2, 5, 6\}& $\emptyset$
\\
\cline{1-2} \cline{4-5} \cline{7-8} \cline{10-11}
                         \{2, 5, 7\}& \{5\}&&
                         \{2, 5, 8\}& \{2\}&&
                         \{2, 5, 9\}& $\emptyset$&&
                         \{2, 6, 7\}& $\emptyset$
\\
\cline{1-2} \cline{4-5} \cline{7-8} \cline{10-11}
                         \{2, 6, 8\}& $\emptyset$&&
                         \{2, 6, 9\}& \{4\}&&
                         \{2, 7, 8\}& $\emptyset$&&
                         \{2, 7, 9\}& \{3\}
\\
\cline{1-2} \cline{4-5} \cline{7-8} \cline{10-11}
                         \{2, 8, 9\}& \{2\}&&
                         \{3, 4, 5\}& \{1\}&&
                         \{3, 4, 6\}& $\emptyset$&&
                         \{3, 4, 7\}& \{1\}
\\
\cline{1-2} \cline{4-5} \cline{7-8} \cline{10-11}
                         \{3, 4, 8\}& $\emptyset$&&
                         \{3, 4, 9\}& \{1\}&&
                         \{3, 5, 6\}& $\emptyset$&&
                         \{3, 5, 7\}& $\emptyset$
\\
\cline{1-2} \cline{4-5} \cline{7-8} \cline{10-11}
                         \{3, 5, 8\}& \{5\}&&
                         \{3, 5, 9\}& $\emptyset$&&
                         \{3, 6, 7\}& \{1\}&&
                         \{3, 6, 8\}& $\emptyset$
\\
\cline{1-2} \cline{4-5} \cline{7-8} \cline{10-11}
                         \{3, 6, 9\}& \{1\}&&
                         \{3, 7, 8\}& \{3\}&&
                         \{3, 7, 9\}& \{4\}&&
                         \{3, 8, 9\}& \{1\}
\\
\cline{1-2} \cline{4-5} \cline{7-8} \cline{10-11}
                         \{4, 5, 6\}& $\{6\}$&&
                         \{4, 5, 7\}& $\emptyset$&&
                         \{4, 5, 8\}& $\emptyset$&&
                         \{4, 5, 9\}& \{5\}
\\
\cline{1-2} \cline{4-5} \cline{7-8} \cline{10-11}
                         \{4, 6, 7\}& $\emptyset$&&
                         \{4, 6, 8\}& $\emptyset$&&
                         \{4, 6, 9\}& $\emptyset$&&
                         \{4, 7, 8\}& $\emptyset$
\\
\cline{1-2} \cline{4-5} \cline{7-8} \cline{10-11}
                         \{4, 7, 9\}& $\emptyset$&&
                         \{4, 8, 9\}& \{4\}&&
                         \{5, 6, 7\}& \{1\}&&
                         \{5, 6, 8\}& \{2\}
\\
\cline{1-2} \cline{4-5} \cline{7-8} \cline{10-11}
                         \{5, 6, 9\}& \{1\}&&
                         \{5, 7, 8\}& $\emptyset$&&
                         \{5, 7, 9\}& \{2\}&&
                         \{5, 8, 9\}& \{1\}
\\
\cline{1-2} \cline{4-5} \cline{7-8} \cline{10-11}
                         \{6, 7, 8\}& $\emptyset$&&
                         \{6, 7, 9\}& $\emptyset$&&
                         \{6, 8, 9\}& \{2\}&&
                       \{7, 8, 9\}& \{1, 3\}
\\
\cline{1-2} \cline{4-5} \cline{7-8} \cline{10-11}
\end{tabular}
\end{align*}

\subsection{Proof of the main result}

We are now in position to prove the first theorem. 

\proof[Proof of Theorem \ref{t:unicity}]  a) Assume that the elements
of the diagonal belong to $\Dr$ then by Remark \ref{r:reduc} and
Section \ref{s:list} we see 
that there exists no $a\in \Z\setminus\{0\}$ such that $\Tr_a$
exists. Therefore by Proposition \ref{p:suite} the solution is unique.
\\
\noindent b) Let $S_0$ be a solution as in \eqref{e:sol}. Suppose now
that the elements of the diagonal does not belong to $\Dr$. Then by
Section \ref{s:list} there is $c>0$  such that
$\Tr_{c}$ and $\Tr_{-c}$ exist. We first discuss the case 
where $c$ is unique. Thanks to \eqref{e:a>0}, we see that we will have
at most one other solution. Suppose now that we have two different
$c>0$ such that $\Tr_{c}$ exists, namely $c_1$ and $c_2$. Lemma
\ref{l:diff} ensures that $\Tr_{c_1}\neq \Tr_{c_2}$. Using again
\eqref{e:a>0}, we conclude we will have
at most one other solution.  \qed 

We conclude by counting the number of unique solutions.

\proof[Proof of Theorem \ref{t:main}] We split the proof according to
the number of $c>0$ that exists in \eqref{e:a+}. 
\\
a) no $c$: there are $35$ different sets $\{x_1, x_2, x_3\}$ that
give that $\Tr_c$ and $\Tr_{-c}$ do not exist (as seen before the solution is
therefore unique). By placing $x_1, x_2,$ and $x_3$ on the diagonal,
we have $3!$ choices and $6!$ choices for the rest of the grid. We
obtain  
\[ 35\times 3!\times 6! \mbox{ possibilities}\]
that lead to a unique solution.
\\
b) A unique $c$: there are $45$ different sets $\{x_1, x_2, x_3\}$ that
give a unique $c$ in \eqref{e:a+}. Then we place $x_1, x_2,$ and $x_3$
on the diagonal. Suppose that we have two solutions (recall that there
can not be more than $2$ by Theorem \ref{t:unicity}). Recalling Remark
\ref{r:sign} we consider solutions either of the form $S_{c}$ or
$S_{-c}$ (but not both of them). Finally using \eqref{e:a>0},  we have
$3! \times 3!$ grids with two solutions.  We obtain  
\[ 45\times 3!\times (6!- 3! \times 3!) \mbox{ possibilities}\]
that lead to a unique solution. 
\\
c) Two different $c$: there are $4$ different sets $\{x_1, x_2, x_3\}$ that
give two different $c$ in \eqref{e:a+}, namely $c_1$ and $c_2$. As
before we start by fixing  $x_1, x_2,$ and $x_3$ on the diagonal. Suppose
that we have two solutions. Then recalling Remark \ref{r:sign} the
solutions are either $\{S_{c_1}, S_{c_2}\}$ or $\{S_{-c_1}, S_{-c_2}\}$
(but not the union). Finally by Lemma \ref{l:diff} and by using 
we have $2\times 3! \times 3!$ grids with two solutions. 
 We
obtain  
\[ 4\times 3!\times (6!- 2\times 3! \times 3!) \mbox{ possibilities}\]
that lead to a unique solution. 
\par
To conclude it remains to sum the three results. \qed

\end{document}